\documentclass[runningheads,a4paper]{llncs}

\usepackage{booktabs}
\usepackage{color}
\usepackage{amsmath,amsfonts,bm,mathtools}
\usepackage{ulem}
\usepackage{amssymb}
\usepackage{graphicx}
\usepackage{algorithmicx}
\usepackage{algpseudocode}

\usepackage{url}
\urldef{\mailsa}\path||
\urldef{\mailsb}\path||
\urldef{\mailsc}\path||
\newcommand{\keywords}[1]{\par\addvspace\baselineskip
\noindent\keywordname\enspace\ignorespaces#1}

\begin{document}

\mainmatter  % start of an individual contribution

\title{A Discrete State Transition Algorithm for Generalized Traveling Salesman Problem}

\titlerunning{Discrete STA for Generalized Traveling Salesman Problem}

\author{Xiaolin Tang \and Chunhua Yang%
\thanks{Corresponding author for this paper. The work is supported by the National Science Found for
Distinguished Young Scholars of China (Grant No. 61025015), Key Project of National Natural Science Funds(Grant No. 61134006), Innovative Research Team in University of Ministry of Education of China(Grant No. IRT1044)}%
\and Xiaojun Zhou\and Weihua Gui}
\authorrunning{Discrete STA for Generalized Traveling Salesman Problem}

\institute{School of Information Science and Engineering, \\ Central South University, Changsha 410083, P.R. China
\mailsa\\
\mailsb\\
\mailsc\\
%\url{http://www.springer.com/lncs}
}


\toctitle{Lecture Notes in Computer Science}
\tocauthor{Authors' Instructions}
\maketitle

\begin{abstract}
Generalized traveling salesman problem (GTSP) is an extension of classical traveling salesman problem (TSP), which is a combinatorial optimization problem and an NP-hard problem. In this paper, an efficient discrete state transition algorithm (DSTA) for GTSP is proposed, where a new local search operator named \textit{K-circle}, directed by neighborhood information in space, has been introduced to DSTA to shrink search space and strengthen search ability. A novel robust update mechanism, restore in probability and risk in probability (Double R-Probability), is used in our work to escape from local minima. The proposed algorithm is tested on a set of GTSP instances. Compared  with other heuristics, experimental results have demonstrated the effectiveness and strong adaptability of DSTA and also show that DSTA has better search ability than its competitors.
\end{abstract}

\section{Introduction}
\indent
Generally speaking, GTSP can be described as follow: given a completely undirected graph $G = \{ V,E\} $, where $V$ is a set of $n$ vertices and has been partitioned into $m$ clusters $V = \left\{ {{V_1},{V_2}, \cdots ,{V_m}} \right\}$, $E$ is a set of $m$ edges, and the goal of GTSP is to find a tour visiting each cluster exactly once while minimizing the sum of the route costs. In this paper, the symmetric GTSP is concerned, that is to say, ${c_{i,j}} = {c_{j,i}}$, here, the associated cost ${c_{i,j}}$ for each pair of vertices $(i,j)$  represents the distance from one vertex in ${V_i}$ to another vertex in ${V_j}$. Since each cluster has at least one vertex and each vertex can only belong to one cluster, we have $m \le n$.  If $m = n$, GTSP is restored to TSP which together with GTSP are both NP-hard problems[1]. To deal with TSP, we only need to optimize the sequence of the clusters, while in the process of solving GTSP we must determine the sequence of the clusters and a vertex to be visited in each cluster simultaneously, which indicates that GTSP is more complex than TSP. Nonetheless, GTSP is extensively used in many applications, such as task scheduling, airport selection and postal routing, etc[2,3].\\
\indent
According to the characteristics of GTSP, we can decompose the process of solving GTSP into two phases. One is to determine the visiting order of all the clusters, which is similar to TSP; the other is to find the optimal vertex in each cluster in a given order. Many reputed heuristic searching algorithms, like genetic algorithm (GA)[4] , particle swarm optimization (PSO)[5], simulated annealing (SA)[6], ant colony optimization (ACO)[7], have been varied into discrete versions to solve GTSP. Though these algorithms have their own mechanisms to deal with continuous optimization problems, they have to adapt themselves to GTSP with some classic operators, such as \textit{swap} and \textit{insert}. These search operators change the visiting order of clusters in particular ways. Lin-Kernighan(L-K) is a well-known method to solve TSP and GTSP[8], which focuses on changing the edges instead of the visiting order of clusters. The number of edges that L-K impacts in a single operation is unknown; as as result, the depth of L-K is usually limited within a constant[9]. The majority of these methods focus on finding an optimal sequence of clusters, while to solve GTSP, it still has to choose a vertex from each cluster to make the minimal cost simultaneously. This is a well-known shortest path problem in operations research which is also called cluster optimization (CO) in GTSP. The most common method to deal with this problem is dynamic programming that can give us a definitively best result which is named as layer network method in other literatures[10].\\
\indent
State transition algorithm is a new optimization algorithm, according to the control theory and state transition[11]. The efficiency of STA in application to continuous optimization problems has been proved[12,13]. In [14,15], Discrete version of state transition algorithm has been introduced to solve a series of discrete optimization problems such as TSP, boolean integer programming. In this study, we will extend DSTA to solve the GTSP.\\
\indent
In section 2, we give a brief description of DSTA and some transformation operators. Section 3 introduces relevancy and correlation index to describe $K$-Neighbor. In section 4, a DSTA is presented to solve GTSP with a new updating mechanism. Some experimental results are given in section 5, and the final part is the conclusion.

\section{Discrete State Transition Algorithm}
\subsection{Description of DSTA}
\indent
State transition algorithm comes from control theory. It regards a solution to an optimization problem as a state and updating of the solution as state transition. The unified form of discrete state transition algorithm is given as follow:
\begin{eqnarray}
\left \{ \begin{array}{ll}
\bm x_{k+1}= A_{k}(\bm x_{k}) \bigoplus B_{k}(\bm u_{k})\\
y_{k+1}= f(\bm x_{k+1})
\end{array} \right.,
\end{eqnarray}
Where, $\bm x_k \in \mathbb{R}^n$ denotes a current state, corresponding to current solution of an optimization problem; ${\bm u_k} \in \mathbb{R}^n$ is a function of ${\bm x_k}$ and historical states; both ${A_k}$, ${B_k} \in {\mathbb{R}^{n \times n}}$ are transition operators which are usually state transition matrixes; $\bigoplus$ is an operation, which is admissible to operate on two states; $f$ is the cost function or evaluation function.\\
\indent
In general, the solution to discrete optimization problem is a sequence, which means a new state ${\bm x_{k + 1}}$ should also be a sequence after transformation by ${A_k}$ or ${B_k}$. For the TSP, only a state transition matrix is considered, avoiding the complexity of “adding” one sequence to another. So the form of DSTA for TSP is simplified as follow:
\begin{equation}
\left\{ \begin{array}{l}
 {\bm x_{k + 1}} = {G_k}{\bm x_k} \\
 {y_{k + 1}} = f({\bm x_{k + 1}}) \\
 \end{array} \right.
\end{equation}
\indent
where, $\bm x_k = {[{x_{1,k}},{x_{2,k}}, \cdots ,x_{m,k}]^T},{x_{i,k}} \in \{ 1,2, \cdots ,m\} $; ${G_k}$ is the state transition matrix which
is created by transformation operators. State transition matrixes are variants of identity matrix with only position value 1 in each column and each row.
Multiplying a state transition matrix by a current state will get a new state which is still a sequence and the process is like this:
%\begin{mdframed}
\begin{multicols}{2}
Current state ${\bm x_k}$: [1 2 3 4 5]$^{\rm{T}}$
\\
\indent
State transition matrix ${G_k}$:
\[\left( {\begin{array}{*{20}{c}}
   1 & 0 & 0 & 0 & 0  \\
   0 & 1 & 0 & 0 & 0  \\
   0 & 0 & 0 & 0 & 1  \\
   0 & 0 & 0 & 1 & 0  \\
   0 & 0 & 1 & 0 & 0  \\
\end{array}} \right)\]
\[\left( {\begin{array}{*{20}{c}}
   1 & 0 & 0 & 0 & 0  \\
   0 & 1 & 0 & 0 & 0  \\
   0 & 0 & 0 & 0 & 1  \\
   0 & 0 & 0 & 1 & 0  \\
   0 & 0 & 1 & 0 & 0  \\
\end{array}} \right)\begin{array}{*{20}{c}}
   {\begin{array}{*{20}{c}}
   {} & {} & {}  \\
\end{array}} &  \times  & {} & {\begin{array}{*{20}{c}}
   {\begin{array}{*{20}{c}}
   {} & {}  \\
\end{array}} & {\left( {\begin{array}{*{20}{c}}
   1  \\
   2  \\
   3  \\
   4  \\
   5  \\
\end{array}} \right)}  \\
\end{array}}  \\
\end{array}\]
New state:[1 2 5 4 3]$^{\rm{T}}\leftarrow{G_k} \times {\bm x_k}$
\end{multicols}
%\end{mdframed}

\subsection{State transformation operators}
\begin{figure}
\centering
\includegraphics[width=12cm,height=8cm]{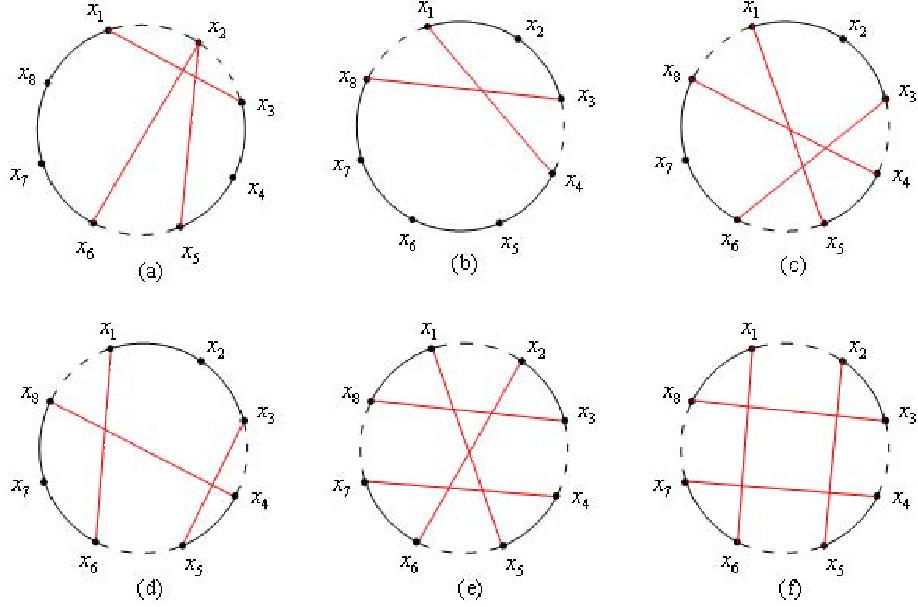}
\caption{Effects of circle operator in different cases}
\label{circles}
\end{figure}

\indent
DSTA for TSP works with 3 efficient operators, which is the foundation of study on GTSP. All of the 3 operators, swap, shift, symmetry, belong to ${G_k}$ .
\\
\noindent
1) Swap
\[({x_1},{x_2},{x_3}, \cdots ,{x_{i - 1}},{x_i},{x_{i + 1}}, \cdots ,{x_{m - 1}},{x_m})\]
\[ \to ({x_1},{x_i},{x_3}, \cdots ,{x_{i - 1}},{x_2},{x_{i + 1}}, \cdots ,{x_{m - 1}},{x_m})\]
\indent
This is an operator to exchange several vertices in the tour and the number of the vertices to be changed is limited by a parameter ${m_a}$. With this operator, the number of edges to be changed is twice as that of vertices to be exchanged.
\\
\noindent
2)	Shift
\[({x_1},{x_2},{x_3}, \cdots ,{x_{i - 1}},{x_i},{x_{i + 1}}, \cdots ,{x_{m - 1}},{x_m})\]
\[ \to ({x_1},{x_3}, \cdots ,{x_{i - 1}},{x_i},{x_2},{x_{i + 1}}, \cdots ,{x_{m - 1}},{x_m})\]
\indent
This operator first removes a segment of sequence from a given tour and then inserts this segment into a random position of the remaining sequence. The length of the removed sequence is restricted to less than ${m_b}$. Three edges will be changed through this operator, of which two edges are adjacent and the last edge is non-adjacent to them.
\\
\noindent
3)	Symmetry
\[({x_1},{x_2}, \cdots ,{x_{i - 3}},{x_{i - 2}},{x_{x - 1}},{x_i},{x_{i + 1}},{x_{x + 2}},{x_{x + 3}} \cdots ,{x_m})\]
\[ \to ({x_1},{x_2}, \cdots ,{x_{i - 3}},{x_{i + 2}},{x_{x + 1}},{x_i},{x_{i - 1}},{x_{x - 2}},{x_{x + 3}} \cdots ,{x_m})\]
\indent
Symmetry is a unique operator in STA, which is to choose a vertex in a given tour as center, and then mirror a small segment on the left side of the center to the opposite side and so does a small segment on the right side. The length of these small segments is restricted to ${m_c}$. For symmetric GTSP, the symmetry operator can change two edges every time.
\\
\noindent
4) Circle\\
\indent
At the final stage of search process, a tour is usually locally optimal, which indicates that this tour has many similar segments to the global tour; therefore, the entire sequence can be regarded as the combination of these segments. To further optimize a tour, we have to change several conjoint vertices in some segments simultaneously.
Considering that the number of changing vertices in swap or shift is no more than 2, a new operator called circle is proposed here to enhance the global search ability.\\
\indent
Circle consists of two steps. First, we divide a given tour into two circles randomly, and then break one of the circles and insert it into another to create a new complete tour. Effects of this operator can be summarized as follow:\\
\indent
1.	One of the circles contains only one vertex(Fig.\ref{circles}.a);\\
\indent
2.	Both of the circles contain more than one vertex, and change the connection at the interfaces of each circle(Fig.\ref{circles}.b);\\
\indent
3.	Both of the circles contain more than one vertex, break one of the circles at its interface and insert it into another from a random position except the interface(Fig.\ref{circles}.c, Fig.\ref{circles}.d).\\
\indent
4.	Both of the circles contain more than one vertex, break one of the circles and insert it into another. Both of the breaking position and inserting position are randomly chosen except two interfaces. Fig.\ref{circles}.e and Fig.\ref{circles}.f show the results of this case.

Obviously, circle is much more flexible than the other operators since it can gain 6 different kinds of cases.
\\
\noindent
5) Cluster optimization(CO)

This is a sole operator to find the best path of given visiting order of clusters. A tour $[{x_{1,k}},{x_{2,k}}, \cdots ,{x_{m,k}}]$ with costs $W({{\bm x_k}})$ will be optimized into a new tour $\bm x_k^{'}$ after running CO, here, $W({\bm x_k^{'}} ) \le W( {{\bm x_k}})$ and $cluster(\bm x_k^{'}) = cluster({\bm x_k})$. In general, few of visiting orders of clusters will be changed from $\bm x_k $ to $\bm x_{k}^{'}$, thus we only need to optimize a small segment around the changed clusters.

\section{K-Neighbor}
To improve the global search ability for large-scale problems in limited computational time, it is necessary to avoid some potential bad search space. In this paper, the correlation index and relevancy are proposed, where correlation index is used to assess the correlation of every two clusters and the relevancy is applied to define K-Neighbor which will guide search direction as heuristic information.\\
\textbf{Definition 1:} Let define the distance between the geometric centers of cluster $i$  and cluster $j$ as ${d_{i,j}}$, the sum of distance from the geometric centers of cluster $i$ to the geometric centers of other clusters as ${d_i}$ and denote ${r_{i,j}}$ as the correlation index of cluster $i$ to cluster $j$:
\begin{equation}
{r_{i,j}} = \frac{{1 - \frac{{{d_{i,j}}}}{{{d_i}}}}}{{n - 1}},
\end{equation}
\[\sum\limits_{j = 1}^n {{r_{i,j}}}  = 1,{d_i} = \sum\limits_{j = 1}^n {{d_{i,j}}}.\]
\\
\textbf{Definition 2:} Given ${r_{i,j}}$ as the correlation index of cluster $i$ to cluster $j$ and ${r_{j,i}}$ as the correlation index of cluster $j$  to cluster $i$, then the relevancy ${p_{i,j}}$ of cluster $i$ to cluster $j$ can be formulated as:
\begin{equation}
{p_{i,j}} = \frac{{{r_{i,j}} \times {r_{j,i}}}}{{\sum\limits_{j = 1}^n {{r_{i,j}} \times {r_{j,i}}} }}.
\end{equation}
Calculating each ${p_{i,j}}$, we can get a relevancy matrix. The $i$th row of the relevancy matrix shows the relevancy of cluster $i$ to other clusters. A big ${p_{i,j}}$ indicates cluster $i$ is with high possibility connecting with cluster $j$. After sorting the relevancy matrix in descending order by row, the top $k$ clusters in row $i$ will be the K-Neighbor of cluster $i$. Using K-Neighbor as heuristic information, the global search ability can be improved significantly.

\section{DSTA for GTSP}
To adapt DSTA to GTSP and to make the algorithm more efficient, K-Neighbor is used as heuristic information to guide the search. Thus, $k$-shift, $k$-symmetry and $k$-circle which are all guided by K-Neighbor are included in DSTA. The core procedure of the DSTA for GTSP can be outlined in pseudocode as follows:
\begin{algorithmic}[1]
\Repeat
    \State {[$Best,Best^{*}$] $\gets$ swap($SE,Best,Best^{*},{m_a}$)}
    \State {[$Best,Best^{*}$] $\gets$ shift($SE,Best,Best^{*},{m_b}$)}
    \State [$Best,Best^{*}$]  $\gets$ k-circle($SE,Best,Best^{*}$,$K$-$Neighbor$)
    \State [$Best,Best^{*}$]  $\gets$ k-symmetry($SE,Best,Best^{*}$,$K$-$Neighbor$)
    \State [$Best,Best^{*}$]  $\gets$ k-shift($SE,Best,Best^{*}$,$K$-$Neighbor$)
\Until{the specified termination criterion is met}
\end{algorithmic}
where, $SE$ is the search enforcement, representing the times of transformation by a certain operator;
$Best$ is the best solution from the candidate state set created by transformation operators; $Best^*$ is the best solution in history. There are five operators in DSTA to optimize the sequence of clusters. A short cluster optimization which only optimizes a small segment (no more than 5 vertices) around the changed vertices is contained in every transformation operator to find a minimum path of a given sequence in further. To escape from local minima, a new robust update mechanism, restore and risk in probability, called double R-Probability for short, is introduced. Risk in probability is to accept a bad solution with a probability $p_1$. To ensure the convergence of DSTA, restore in probability $p_2$ is designed to recover the best solution in history.
\section{Computational Results}
Instances used in this paper all come from GTSPLIB[16]. The number of clusters in these instances varies from 30 to 89.  All of the algorithms, including DSTA, SA and ACO, are coded in matlab and run on an Intel Core i5 3.10GHz under Window XP environment.
In order to test the performance of the proposed operators and approach, DSTA is compared with SA and ACO, and 10 runs are carried out for the experiment. Some statistics are used as follows

Opt. : the best known solution,\\
\indent
Best:  the best solution obtained from the experiment,\\
\indent
${\Delta _{avg}}$: the relative error of the average solution,\\
\[{\Delta _{avg}} = \frac{{mean(values) - Opt.}}{{Opt.}} \times 100\% \]
\indent
${t_{avg}}$: average time consumed.

Results of comparison among DSTA, SA and ACO are listed in Table \ref{tab1}. In DSTA, we set $k=8$, $m_a = 2, m_b = 1$. The initial temperature in SA is 5000 and the cooling rate is 0.97. In ACO, $\alpha  = 1$, $\beta  = 5$, $\rho  = 0.95$, where, $\alpha$, $\beta$ are used to control the relative weight of pheromone trail and heuristic value, and $\rho$ is the pheromone trail decay coefficient. As can be seen from Table \ref{tab1}, DSTA is superior to SA and ACO in both time consumption and solution quality. The ${\Delta _{avg}}$ of DSTA is very small, which indicates DSTA has good robustness and can obtain good solutions with high probability. In the 10 runs, DSTA obtains the optimal solution at almost each run for every instance, but SA and ACO seldom find the Opt. except for 30kroB150. ${\Delta _{avg}}$ of SA is smaller than that of ACO because SA accepts a bad solution with probability which can help it escape from local minima.
\begin{table}[h]
\centering  % 表居中
\caption{ Results of comparison for SA, ACO and DSTA.}
\label{tab1}
\begin{tabular}{{p{1.8cm}p{1cm}p{1cm}p{1cm}p{0.8cm}p{1cm}p{1cm}p{0.8cm}p{1cm}p{1cm}p{0.8cm}}}\toprule[0.75px]  % {lccc} 表示各列元素对
\multicolumn{2}{c}{} & \multicolumn{3}{c}{SA} & \multicolumn{3}{c}{ACO}& \multicolumn{3}{c}{DSTA}\\\cmidrule(r){3-5}\cmidrule(r){6-8}\cmidrule(r){9-11}
Instance & Opt.  & Best & ${\Delta _{avg}}$  & $t_{avg}$ & Best & ${\Delta _{avg}}$ & $t_{avg}$ & Best & ${\Delta _{avg}}$ & $t_{avg}$\\ \hline
30kroA150	&	11018	&	11027 	&	0.16 	&	152 	&	11331	&	5.99 	&	104 	&	11018	&	0	&	13 	\\
30kroB150	&	12196	&	12196 	&	0.02 	&	78 	&	12532	&	6.02 	&	67 	&	12196	&	0.18 	&	18 	\\
31pr152	&	51576	&	51584 	&	1.12 	&	79 	&	51734	&	1.60 	&	69 	&	51576	&	0 	&	25 	\\
32u159	&	22664	&	22916 	&	1.90 	&	89 	&	24285	&	8.68 	&	75 	&	22664	&	0.74 	&	33 	\\
39rat195	&	854	&	857 	&	1.09 	&	198 	&	884	&	5.86 	&	145 	&	854	&	0.05 	&	56 	\\
40d198	&	10557	&	10574 	&	0.53 	&	112 	&	11458	&	9.31 	&	103 	&	10557	&	0.06 	&	74 	\\
40kroA200	&	13406	&	13454 	&	0.62 	&	107 	&	14687	&	10.77 	&	99 	&	13406	&	1.10 	&	69 	\\
40kroB200	&	13111	&	13117 	&	0.38 	&	108 	&	13396	&	8.34 	&	99 	&	13111	&	0.20 	&	28 	\\
45ts225	&	68340	&	68401 	&	1.57 	&	325 	&	70961	&	5.83 	&	223 	&	68340	&	0.66 	&	72 	\\
45tsp225	&	1612	&	1618 	&	1.77 	&	122 	&	1736	&	8.15 	&	119 	&	1612	&	1.35 	&	88 	 \\
46pr226	&	64007	&	64062 	&	2.70 	&	130 	&	66458	&	7.51 	&	124 	&	64007	&	0 	&	35 	\\
53gil262	&	1013	&	1047 	&	5.24 	&	142 	&	1148	&	15.92 	&	148 	&	1013	&	1.30 	&	85 	 \\
53pr264	&	29549	&	29725 	&	1.87 	&	146 	&	32388	&	12.06 	&	150 	&	29546	&	0.07 	&	58 	\\
60pr299	&	22615	&	23186 	&	7.00 	&	165 	&	25296	&	15.97 	&	184 	&	22618	&	2.54 	&	114 	 \\
64lin318	&	20765	&	21528 	&	5.73 	&	166 	&	23365	&	13.57 	&	199 	&	20769	&	2.62 	&	117 	 \\
80rd400	&	6361	&	6920 	&	10.36 	&	225 	&	8036	&	21.67 	&	299 	&	6361	&	2.52 	&	141 	 \\
84fl417	&	9651	&	10099 	&	9.95 	&	282 	&	10122	&	10.14 	&	345 	&	9651	&	0.51 	&	158 	 \\
88pr439	&	60099	&	66480 	&	13.13 	&	276 	&	69271	&	16.14 	&	368 	&	60099	&	2.95 	&	156 	 \\
89pcb442	&	21657	&	23811 	&	11.15 	&	253 	&	26233	&	19.48 	&	376 	&	21664	&	3.80 	&	163 	 \\
\bottomrule[0.75px]
\end{tabular}
\end{table}

\section{Conclusions}
We added a new operator and heuristic information to DSTA to solve GTSP. K-Neighbor can guide the search direction, in a way to ignore all possible connections among vertices.
A flexible operator $k$-circle is guided by the K-Neighbor, which can change random segments freely in a tour. Double R-Possibility is helpful to escape from local minima. It accepts a bad solution with a probability ${p_1}$ and restore the history best with another probability ${p_2}$. All these strategies contribute
to improving the performance of the DSTA.

\label{references}

\end{document}